\documentclass[11pt]{article}
\usepackage{harvard}

\usepackage{amssymb}
\usepackage{amsmath}
\usepackage{amsthm}
\usepackage{mathrsfs}
\usepackage{comment}
\usepackage{a4}

\usepackage{graphics}
\usepackage{float}
\usepackage{epsfig}

\usepackage{dsfont}
\usepackage{epic}
\usepackage{eepic}
\usepackage{pgf}

\newtheorem{theorem}{Theorem}[section]
\newtheorem{lemma}[theorem]{Lemma}
\newtheorem{proposition}[theorem]{Proposition}

\newtheorem{definition}[theorem]{Definition}
\newtheorem{assumption}[theorem]{Assumption}
\newtheorem{remark}[theorem]{Remark}
\newtheorem{example}[theorem]{Example}

\renewcommand{\(}{\left(}
\renewcommand{\)}{\right)}

\newcommand{\argmin}{\operatorname*{argmin}}

\DeclareMathOperator{\E}{{\mathbb E}}

\DeclareMathOperator{\R}{{\mathds R}}

\DeclareMathOperator{\N}{{\mathds N}}

\DeclareMathOperator{\PP}{{\mathbb P}}

\providecommand{\scapro}[2]{\langle #1,#2 \rangle}
\providecommand{\eps}{\varepsilon} \providecommand{\norm}[1]{\lVert
#1 \rVert} \providecommand{\abs}[1]{\lvert #1 \rvert}

\allowdisplaybreaks[2]

\providecommand{\babs}[1]{{\Bigl\lvert #1 \Bigr\rvert}}

\begin{document}

\title{Regularization independent of the noise level:\\ an analysis of quasi-optimality}
\author{\parbox{5.5cm}{\centering Frank Bauer\\Fuzzy Logic Laboratory\\University of
Linz\\ frank.bauer@jku.at} \hspace{1mm}
\parbox{6.5cm}{\centering Markus Rei{\ss}\\ Institute of Applied
Mathematics\\University of Heidelberg\\
reiss@statlab.uni-heidelberg.de} }

\maketitle


\begin{abstract}
The quasi-optimality criterion chooses the regularization parameter in inverse problems
without taking into account the noise level. This rule works remarkably well in
practice, although Bakushinskii has shown that there are always counterexamples with
very poor performance. We propose an average case analysis of quasi-optimality for
spectral cut-off estimators and we prove that the quasi-optimality criterion determines
estimators which are rate-optimal {\em on average}. Its practical performance is
illustrated with a calibration problem from mathematical finance.
\end{abstract}

\section{Introduction}

We consider the prototype of a linear inverse problem where we observe $y=Ax+\xi$ with a
compact operator $A$ on a Hilbert space $X$ and some noise variable $\xi$ and where we
try to recover a stable approximation of the true solution $x\in X$. In this setting we
address the question, why the so-called quasi-optimality criterion for choosing the
regularization parameter in stable inversion algorithms, as proposed by
\citeasnoun{Tikhonov/Arsenin:1977} and \citeasnoun{Tikhonov/Glasko/Kriksin:1979}, works
remarkably well in many practical situations.

The classical quasi-optimality criterion is applied to the
regularized solutions
\begin{equation*}
x_\alpha = \( A^* A + \alpha I \)^{-1} A^* y,\quad \alpha>0,
\end{equation*}
based on Tikhonov's method and proposes to choose $\alpha>0$ such
that
\begin{equation*}
\left\|\alpha \frac{d x_\alpha}{d \alpha} \right\| \rightarrow
\min_\alpha!
\end{equation*}
Reparametrizing with $\alpha=q^n$ for $q\in(0,1)$ and $n\in\R$
yields $\| \frac{d x_{q^n}}{d n}\| \rightarrow \min_n!$ When for
practical purposes only a discrete grid $\{q^n\,|\,n\in\N\}$ is
considered, then this reduces to the following criterion
\begin{equation*}
\| x_{q^n} - x_{q^{n+1}} \|\rightarrow \min_n!
\end{equation*}
Obviously this parameter choice rule does not rely on any knowledge
concerning the operator, the solution and the noise variable.

On the other hand, \citeasnoun{Bakushinskii:1984} has shown for deterministic noise
$\xi$ that a method for choosing the regularization parameter should depend on the noise
level, whereas quasi-optimality does not.

\begin{theorem}
If the regularized solution operator $R:X\to X$ does not depend
explicitly on the noise level $\delta$, then for any (infinite rank)
compact operator $A:X\to X$ there is a $y \in \operatorname{Dom}
(A^+)$ for which
\begin{equation*}
\lim_{\delta \rightarrow 0} \sup_{\begin{smallmatrix} y_\delta\in X
\\ || y -y_\delta || \leq \delta \end{smallmatrix}} || R y_\delta -
A^+ y ||
> 0,
\end{equation*}
where $A^+$ denotes the generalized inverse of $A$.
\end{theorem}
The core message of our paper is that while Bakushinskii's result is
true for a worst case analysis, that is for any such method $R$
certain counterexamples can be constructed, it is not necessarily
true when we assess a method by its average case performance. In
particular, quasi-optimality works well on average. Our analysis
applies also to variants of quasi-optimality like Hardened Balancing
\cite{Bauer:2007}, which are sometimes preferable in practice.

For a general overview about practical inverse problems, their abstract mathematical
formulation and methods for choosing the regularization parameter we refer to the
monograph \citeasnoun{Engl/Hanke/Neubauer:1996}, while the Bayesian framework, which is
related to our approach, is discussed by \citeasnoun{Kaipio/Somersalo:2005}. The choice
of the regularization parameter is, of course, a perennial problem in nonparametric
statistics, see e.g. \citeasnoun{Wasserman:2006}, and diverse methods have been applied
to statistical inverse problems, e.g. generalized cross validation \cite{Wahba:1977},
Lepski's method \cite{Lepskij:1990,Bauer/Pereverzev:2005} or wavelet thresholding
\cite{Cohen/Hoffmann/Reiss:2004}, which depend, however, heavily on the knowledge of the
noise level.  The easy to use and heuristically motivated quasi-optimality criterion and
its variants are attractive in practice because they do not rely on the knowledge of the
noise level.

First considerations, why this kind of methods might work, are already given by
\citeasnoun{Bauer:2007}. In Section \ref{SectionMathAn} below, we derive the proper
mathematical result that estimators, based on a spectral cut-off scheme and the
quasi-optimality criterion for selecting the cut-off value, are rate-optimal {\em on
average}. Specifying an average case scenario amounts to prescribing a Bayesian a priori
law for the functions to be estimated. We suppose that the coefficients of the solution
in the singular value decomposition are normally distributed around zero. Only a very
general condition on the decay property of the variances will be imposed for the proof
of the theorem. The precise setting is described in Section \ref{SectionFramework},
which is followed by the mathematical analysis in Section \ref{SectionMathAn}. The
proofs are postponed to Section \ref{SectionProofs}. The numerical example in Section
\ref{SectionAppl} treats the calibration an option price model from mathematical
finance. In this typical inverse problem in financial engineering the choice of the
regularization parameter is extremely difficult and the proposed methods yield
comparatively good results. A short outlook in Section \ref{SectionPersp} concludes.

\section{The framework}\label{SectionFramework}

\subsection{Observations and estimators}

We consider $A:X\to X$, a compact, self-adjoint and
positive-definite operator on the real Hilbert space $X$ with
singular value decomposition
\begin{equation}
Ax=\sum_{k=1}^\infty \lambda(k)\scapro{x}{u_k}u_k =
\sum_{k=1}^\infty \lambda(k)x_ku_k,
\end{equation}
where $(u_k)$ is an orthonormal basis of eigenvectors and the
positive eigenvalues $(\lambda(k))_k$ are arranged in decreasing
order, satisfying $\lim_{k\to\infty}\lambda(k)=0$. From the general
observation model
\[ y=Ax+\xi,\quad\text{$\xi$ noise},\]
we immediately go over to a sequence space model by considering the
coordinates with respect to $(u_k)$. We observe
\[ y_k=\lambda(k)x_k+\eps(k)\xi_k,\quad k\ge 1.\]
Here, $\xi_k$ are i.i.d. standard normal random variables and
$(x_k)$ the coordinates of the unknown quantity $x$, which is to be
estimated. Writing $\sigma(k):=\eps(k)/\lambda(k)$, the empirical
coefficients of $x$ are given by
\begin{equation}\label{EqModel}
\tilde
x_k:=\lambda(k)^{-1}y_k=x_k+\eps(k)\lambda(k)^{-1}\xi_k=x_k+\sigma(k)\xi_k,\quad
k\ge 1.
\end{equation}
Equation \eqref{EqModel} is the abstract observation model we shall
consider from now on. We shall use a subsampling function
$\ell:\N\to\N$ with $\ell(n+1)>\ell(n)$. Applying a spectral cut-off
scheme, our estimator of $x$ at level $n$ for a given subsampling
$\ell$ is defined as
\[ \hat{x}^{(n)}:=\sum_{k=1}^{\ell(n)}\tilde x_k u_k.\]
Its mean squared error (MSE) is given by the bias-variance
decomposition, see e.g. \citeasnoun{Wasserman:2006},
\citeasnoun{Engl/Hanke/Neubauer:1996}:
\[ \E[\norm{\hat{x}^{(n)}-x}^2]=\sum_{k=1}^{\ell(n)}
\sigma(k)^2+\sum_{k=\ell(n)+1}^\infty x_k^2.
\]
Recall that the bias-variance dilemma is the fact that the variance,
i.e. the first sum, increases with $n$, while the squared bias, i.e.
the second sum, decreases. The value of $n$ where the total sum is
minimal depends on bias and variance and thus on the properties of
the unknown $x$ and of the noise level $\sigma(\cdot)$. The
quasi-optimality criterion gives a data-driven choice for $n$.

\subsection{Assumptions}

Let us now adopt a Bayesian point of view and perform an average
case analysis. We weight the coefficients of $x$ by the prior
distribution
\[ x_k\stackrel{i.i.d.}{\sim} N(0,\gamma(k)^2),\quad k\ge 1.\]
Assumption \ref{AssDecay} below will implicitly require certain
decay properties of the variances $\gamma(k)^2$ for $k\to\infty$,
but we need not know them in detail. Writing $\tilde\E$ for the
joint expectation with respect to $(\xi(k))$ and $(x_k)$, the
Bayesian risk $R_2$ for the MSE is given by
\[ R_2(\hat{x}^{(n)})^2:=\tilde\E[\norm{\hat{x}^{(n)}-x}^2]
=\sum_{k=1}^{\ell(n)} \sigma(k)^2+\sum_{k=\ell(n)+1}^\infty
\gamma(k)^2.
\]
Introducing
\begin{align*}
s(n)&:=\sum_{k=1}^{\ell(n)} \sigma(k)^2 &&\text{(variance)},\\
b(n)&:=\sum_{k=\ell(n)+1}^\infty \gamma(k)^2 &&\text{(mean squared
bias)}
\end{align*}
yields the risk decomposition $R_2(\hat{x}^{(n)})^2=s(n)+b(n)$.

\begin{assumption}\label{AssDecay} (Geometric growth, decay rates)
We assume that with some $1<c_s\le C_s$, $1<c_b\le C_b$ we have for
all $n\ge 1$
\[ c_s s(n)\le s(n+1)\le C_s s(n),\quad
c_b b(n+1)\le b(n)\le C_b b(n+1).
\]
\end{assumption}
This assumption can be easily fulfilled for most moderately
ill-posed problems with H\"older source conditions
\cite{Engl/Hanke/Neubauer:1996} by choosing an exponential
subsampling function $\ell$.

\begin{example}\label{RemarkExample}
Assume that we have a moderately ill-posed problem with H\"older
source conditions (i.e. $\lambda(k) \asymp k^{-\nu}$, $\nu>0$, and
$\gamma(k) \asymp k^{-\mu}$, $\mu
>1/2$) and white noise of level $\delta>0$ (i.e. $\eps(k) = \delta$, and hence $\sigma(k) = \delta
\lambda(k)^{-1} \asymp \delta k^{\nu}$). This white noise setting is
assumed for simplicity, but is not strictly necessary.

Then choosing an exponential subsampling like  $\ell(n)=\ell(0) h^n$
with $\ell(0)\in\N$, $h>1$ independent of the noise level and the
source conditions, we obtain
\begin{align*}
s(n)&=\sum_{k=1}^{\ell(n)} \sigma(k)^2\asymp\sum_{k=1}^{\ell(0) h^n}
\delta^2 k^{2\nu} \asymp
\tfrac{\delta^2\ell(0)^{2\nu+1}}{2\nu+1}  (h^{2\nu+1})^n  \\
b(n)&=\sum_{k=\ell(n)+1}^\infty \gamma(k)^2 \asymp\sum_{k=\ell(0)
h^n+1}^\infty k^{-2\mu} \asymp \tfrac{\ell(0)^{-2\mu+1}}{2\mu -1}
(h^{-2\mu+1})^n
\end{align*}
and hence $s(n)$ and $b(n)$ fulfill Assumption \ref{AssDecay}.
\end{example}

Right now we will introduce a weight function $\chi$, which will
allow to treat also variants of quasi-optimality, cf. Remark
\ref{remark} below.

\begin{assumption}\label{AssDecay2} (weight function)
Let $\chi:\N\to\R$ be a weight function which satisfies for some
constants $c_\chi,C_\chi$ with $c_b^{-1}<c_\chi\le C_\chi<c_s$:
\begin{equation*}
c_\chi \chi(n+1) \le \chi(n) \le C_\chi \chi(n+1).
\end{equation*}
\end{assumption}

\begin{example}
For $\chi(\cdot) = 1$ Assumption \ref{AssDecay2} is obviously
fulfilled because of Assumption \ref{AssDecay}. Using that $s(n)$ is
increasing, Assumption \ref{AssDecay2} also holds for $\chi(n) =
s(n)^{-1/2}$.
\end{example}

\subsection{Choosing the regularization parameter}

\begin{definition}\label{DefParameter}
Given a weight function $\chi$, we choose the cut-off level in a
data-driven way according to the minimum distance or
quasi-optimality criterion
\cite{Tikhonov/Arsenin:1977,Tikhonov/Glasko/Kriksin:1979}:
\[ n^\ast:=\argmin_{n\ge 1} \Big\{\chi(n) \norm{\hat{x}^{(n+1)}-\hat{x}^{(n)}}^2\Big\}.\]
\end{definition}

\begin{remark}\label{remark}
For $\chi(\cdot) = 1$ this is a discretized version of the
quasi-optimality criterion; for $\chi(n) = s(n)^{-1/2}$ a version of
the hardened balancing principle \cite{Bauer:2007}. In Theorem
\ref{MainTheorem} it is proved that the infimum of the criterion
over $n$ is almost surely attained and thus $n^\ast$ is well
defined. It is unique when we take as minimizer the smallest index
$n$ where the minimum is attained.
\end{remark}

Obviously, we do not need at any point an explicit or implicit
knowledge of the noise level for computing
$\norm{\hat{x}^{(n+1)}-\hat{x}^{(n)}}$.

The question remains how to minimize the criterion numerically.
First of all, we have in practice an idea about a lower bound for
the noise (e.g. the machine precision). Furthermore, in practical
applications we only have a finite number of observations and we
shall never deal with more coefficients in the sequence space model
than observations available.

\subsection{Intuition for the proof}

\begin{figure}
\pgfdeclareimage[interpolate=false,width=\linewidth]{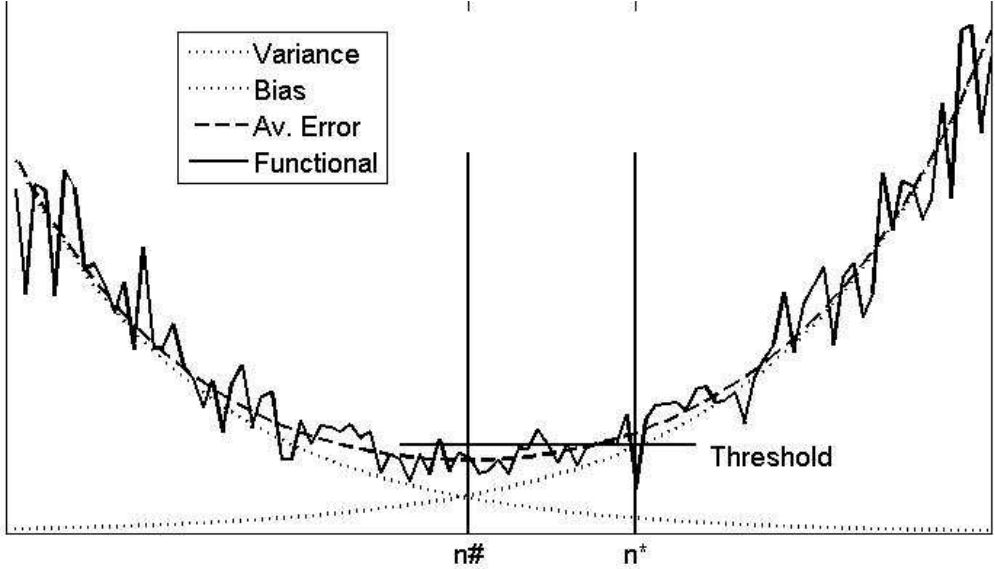}{example}
\pgfuseimage{example}\caption{Idea of Proof} \label{fig0}
\end{figure}

Let us consider the case $\chi(\cdot)=1$. At a first heuristic level
the quasi-optimality criterion is plausible because by the geometric
decay and growth in Assumption \ref{AssDecay}:
\[
\tilde\E[\norm{\hat{x}^{(n)}-\hat{x}^{(n+1)}}^2]=(s(n+1)-s(n))+(b(n)-b(n+1))\asymp
s(n)+b(n)=R_2(\hat{x}^{(n)})^2.
\]
Hence, minimizing the criterion is related to minimizing the average
case risk. Note, however, that a careful analysis is needed since
the criterion is random and not at all independent of the estimator.

More generally, Definition \ref{DefParameter} can be understood as a search for the
intersection point $n^\#$ of the decreasing function $b(\cdot)$ and the increasing
function $s(\cdot)$, where the minimal risk (with respect to $n$) is attained up to some
constant factor, see Figure \ref{fig0} and Lemma \ref{LemIneq} below. Usually,
comparably much is known about the variance part $s(\cdot)$ and not so much about the
bias part $b(\cdot)$.

Let us briefly explain in words why the data-driven choice $n^\ast$
will be close to the intersection point $n^\#$ and therefore the
resulting error in estimating $x$ does not change in order. As we
face a symmetric situation (by exchanging the role of $s(\cdot)$ and
$b(\cdot)$), we consider $n^\ast$ on the right side of $n^\#$
($n^\ast>n^\#$).

In expectation (the situation depicted in Figure \ref{fig0}) the
error curves meet at the oracle value $n^\#$. Here, though, we need
to compute the probability that another point $n^*$ yields the
minimal point for $\norm{\hat{x}^{(n+1)}-\hat{x}^{(n)}}$.

To bound this probability, we introduce a threshold; the probability
that $\norm{\hat{x}^{(n^\#+1)}-\hat{x}^{(n^\#)}} >
\text{threshold}(n^*)$ or $\norm{\hat{x}^{(n^* +1)}-\hat{x}^{(n^*)}}
< \text{threshold}(n^*)$ is larger than the desired probability for
$\norm{\hat{x}^{(n^\#+1)}-\hat{x}^{(n^\#)}} >
\norm{\hat{x}^{(n^*+1)}-\hat{x}^{(n^*)}}$ but much easier
computable.

Furthermore, we completely ignore the fact we can just have one
minimum, we replace $\tilde\E[\norm{x-\hat{x}^{(n)}}^2]$ by the sum
of all positions where $n^*$ beats $n^\#$ weighted by the
probability of their occurrence. Even for this very rough estimation
this sum is still bounded from above by the order of the oracle. The
formal mathematical statements are formulated along these lines in
the next section.

\section{The mathematical analysis}\label{SectionMathAn}

In this section we first determine a rate-optimal estimator for
different moment-type losses, based on the (unrealistic) knowledge
of the index $n^\#$ where the curves $b(\cdot)$ and $s(\cdot)$
intersect. Then we show that this rate does not deteriorate when we
take the estimator with the data-driven index $n^\ast$, based on the
quasi-optimality criterion of Definition \ref{DefParameter}. For the
sake of readability all proofs are deferred to Section
\ref{SecProofProp}.

\subsection{Risk for different moments}

We denote the normalized Bayes risk for the moments of order
$\alpha>0$ by
\[
R_{\alpha}(\hat{x}^{(n)}):=\E[\norm{\hat{x}^{(n)}-x}^{\alpha}]^{1/\alpha},
\]
consistent with the definition of $R_2(\hat{x}^{(n)})$. Note the
order $R_{\alpha}(\hat{x}^{(n)})\le R_{\beta}(\hat{x}^{(n)})$ for
$\alpha\le\beta$.

\begin{definition}
Let $n^\#$ be the index where $s(n^\#+1)> b(n^\#+1)$ and $s(n^\#)\le
b(n^\#)$.
\end{definition}
Note that $n^\#$ is uniquely defined by the monotonicity properties
of $s(\cdot)$ and $b(\cdot)$. Interestingly, for $n^\#$ the risk is
up to a multiplicative factor minimal, even under different moments:

\begin{lemma}\label{LemIneq}
Denote by $\Gamma$ the Gamma function and by $\zeta$ a standard
normal random variable. With
\[ K_{\alpha}:= \sqrt{2 C_b}\,(4\Gamma(\tfrac{\alpha}{2}+1))^{1/\alpha} e^{-\E[\log(\abs{\zeta})]}
\]
the following bounds hold for all $\alpha>0$:
\begin{align*}
\forall n\in\N:\; R_{\alpha}(\hat{x}^{(n)})&\ge
e^{\E[\log(\abs{\zeta})]}R_2(\hat{x}^{(n)}),\\
\forall n\in\N:\; R_{\alpha}(\hat{x}^{(n)})&\le (4
\Gamma(\tfrac{\alpha}{2}+1))^{1/\alpha}R_2(\hat{x}^{(n)}),\\
\min_nR_{\alpha}(\hat{x}^{(n)})&\le R_{\alpha}(\hat{x}^{(n^\#)}) \le
K_{\alpha}\min_nR_{\alpha}(\hat{x}^{(n)}).
\end{align*}
\end{lemma}
We conclude from the last property that the risk (the expectation of
the error) at the intersection point of $s(\cdot)$ and $b(\cdot)$ is
bounded by a constant multiple of the minimal possible risk. The
quasi-optimality criterion aims at recovering this index $n^\#$.

\subsection{Main result}

The first proposition quantifies the probability
$\tilde\PP(n^\ast=n)$.

\begin{proposition}\label{PropLD}
Grant Assumptions \ref{AssDecay} and \ref{AssDecay2}. Then
\[ \tilde\PP(n^\ast=n)\le
(\sqrt{2}+(2er(n))^{r(n)/2})(\rho(n)\log(\rho(n)^{-1}))^{r(n)/2},
\]
where
\[r(n):=\frac{\abs{\chi(n)\left((s(n+1)-s(n))+(b(n)-b(n+1))\right)}}{\max_{\ell(n)<k\le\ell(n+1)}(\sigma(k)^2+\gamma(k)^2)}
\]
and
\[\rho(n):=
\begin{cases}
\frac{C_bC_s-1}{c_s-1} (c_s C_\chi^{-1})^{-\abs{n-n^\#}}, &\text{for
}n\ge
n^\#,\\
\frac{C_bC_s-1}{c_b-1} (c_b c_\chi)^{-\abs{n-n^\#}}, &\text{for
}n<n^\#.
\end{cases}
\]
\end{proposition}

Since $\rho(n)$ decays exponentially fast, this result yields a
rapid decay for the probability that the minimal value for the
quasi-optimality criterion is obtained away from the index $n^\#$.
The slope of the exponential is largely depending on the subsampling
function used, which will be investigated further in the upcoming
example. We are now prepared to state and prove our main theorem.

\begin{theorem}\label{MainTheorem}
Grant Assumptions \ref{AssDecay} and \ref{AssDecay2} and set
$r:=\inf_nr(n)$. Assume that $\alpha>0$ satisfies
\[
\alpha<r\min\Big(\frac{\log(c_b c_\chi)}{\log(C_b)}, \frac{\log(c_s
/ C_\chi)}{\log(C_s)}\Big).
\]
Then our estimator is almost surely well defined, i.e. the minimum
in Definition \ref{DefParameter} for $n^\ast$ is indeed attained,
and it satisfies the oracle-type inequality
\[ R_{\alpha}(\hat{x}^{(n^\ast)})\le K
\min_n R_{\alpha}(\hat{x}^{(n)})
\]
with a constant $K=K(c_s,c_b,C_s,C_b,c_\chi,C_\chi,r,\alpha)>0$.
\end{theorem}

\begin{remark}
A closer look at the proof of Proposition \ref{PropLD} reveals that
for this result only the boundedness and the decay behaviour of the
Gaussian density was used. Hence, a similar result will hold for
more general noise densities with exponential decay.
\end{remark}

This oracle-type result, see e.g. \citeasnoun{Wasserman:2006} for a
general discussion, is actually better than the classical statement
obtained in inverse problems. There the results are given in the
form \cite{Engl/Hanke/Neubauer:1996}
\begin{equation*}
\| x - \hat{x}^{(n^\ast)} \| \leq c \delta^\tau,
\end{equation*}
where $\tau$ is dependent on the type of noise and the source
conditions. When $x$ or the noise have ``accidentally'' better
properties than expected we still have this bound although $\min_n
\| x - \hat{x}^{(n)} \| \ll c \delta^\tau$ for $\delta < \delta_0$.

In our case of an oracle inequality this cannot happen because we
always compare with the best possible (``oracle'') solution. So we
still have the theoretical upper bound in the form $c \delta^\tau$
but can additionally claim, that (in average) we are very close to
the best possible solution, even if it is much better than any
specific rate derived.

Although the mathematical result depends on a number of constants,
it is important to note that the existence of these constants is
required to derive the mathematical result, their actual values are
not at all used for the quasi-optimality criterion.

It is interesting to note that though we face a stochastic noise we
do not lose a logarithmic factor in the data-driven method as is
sometimes the case for other parameter choice methods
\cite{Cohen/Hoffmann/Reiss:2004,Bauer/Pereverzev:2005}.

The general picture is that our data-driven choice $n^\ast$ of the
cut-off value yields an optimal risk bound up to constants whenever
the moments taken are not so large. In order to achieve higher
moments, we need to choose a subsampling function $\ell(n)$ that
makes $r$ large enough. In any case, we have $r\geq 1$ by
definition.

\begin{example}
For the standard case of quadratic risk we need $\alpha=2$. For very
tight exponential bounds $C_s\approx c_s$, $C_b\approx c_b$ and
$\chi(\cdot)=1$ it suffices to have $r>2$. Now we can reconsider
Example \ref{RemarkExample}, where
\begin{align*}
a(n) \approx \frac{\delta^2\ell(0)^{2\nu+1}}{2\nu+1}  (h^{2\nu+1})^n
(h^{2\nu+1} -1) + \frac{\ell(0)^{-(2\mu-1)}}{2\mu -1}
(h^{-2\mu+1})^{n+1} (h^{2\mu-1} -1 )
\end{align*}
Hence we have
\begin{align*}
r(n) \approx & \min \left\{
\frac{\frac{\delta^2\ell(0)^{2\nu+1}}{2\nu+1} (h^{2\nu+1})^n
(h^{2\nu+1} -1)}{\delta^2 (\ell(0) h^{n+1})^{2\nu}} ,
\frac{\frac{\ell(0)^{-(2\mu-1)}}{2\mu -1} (h^{-2\mu+1})^{n+1}
(h^{2\mu-1} -1 )}{ (\ell(0) h^n)^{-2\mu} }
\right\} \\
= & \ell(0) \min \left\{ \frac{1}{2\nu+1} h^{n-2\nu} (h^{2\nu+1} -1)
, \frac{1}{2\mu -1} h^{n -2\mu+1} (h^{2\mu-1} -1 ) \right\}
\end{align*}
and so in the worst case
\begin{equation*}
2 < r \leq r(1) \approx \ell(0) h (h-1)  = \ell(2) - \ell(1).
\end{equation*}
Thus the first (and smallest) distance $\ell(2) - \ell(1)$ needs to
be at least $3$ in this case. This finding corresponds very well
with numerical experience.
\end{example}


\section{Application}\label{SectionAppl}

In \citeasnoun{Bauer:2007} the quasi-optimality, Lepski balancing
and hardened balancing principles have been compared numerically.
The findings are that both, for a large scale stochastic experiment
inverting ill-conditioned matrices and for a more realistic
experiment determining the field of gravity from satellite data,
quasi-optimality and hardened balancing perform quite well and
stably, in particular better than the Lepski balancing principle.

The advantage of quasi-optimality is that it can cope with a rather
unknown structure of the noise. Therefore we will present in the
sequel numerical experiments for an inverse problem arising in
option pricing where noise enters from various sources and its level
is not easy to estimate.

The calibration of financial models based on option prices has
attracted increasing attention recently due to its practical
importance and mathematical challenges, see e.g.
\citeasnoun{Crepey:2003} and \citeasnoun{Egger/Hein/Hofmann:2006}
and the references therein for the case of a generalized
Black-Scholes model and Chapter 13 in \citeasnoun{Cont/Tankov:2004}
for jump models as considered here.

\subsection{The calibration problem}

We consider the problem of calibrating an exponential L\'evy model
based on market prices of European options and closely follow
\citeasnoun{Belomestny/Reiss:2006}. We briefly describe the problem,
but refer to \citeasnoun{Cont/Tankov:2004} for a thorough
introduction. It is assumed that we observe prices $C(K_j)$ of
European call options with different strike prices $K_j$,
$j=1,\ldots,n$, and same time $T>0$ to maturity and that the
underlying stock price follows an exponential L\'evy process
\[ S_t=S_0 e^{L_t}\text{ with a L\'evy process $(L_t)$.}\]
Here, the L\'evy process is restricted to be a superposition of a
Brownian motion of volatility $\sigma^2>0$, a linear drift of slope
$\gamma\in\R$ and a compound Poisson jump process of intensity
$\lambda>0$ with jump density $\nu:\R\to\R^+$. The goal is to
estimate these model parameters, in particular the jump density,
which because of only finitely many observations and the presence of
bid-ask spreads is a typical inverse problem with noisy observations
in quantitative finance. The knowledge of these parameters permits
to get a clear picture of the expectations at the market concerning
future jumps in the stock price, which is essential for well-founded
risk management and pricing of path-dependent options.

We transform the strike prices $(K_j)$ to the so called {\em
log-forward moneyness} $(x_j)$ and the call option prices $C(K_j)$
to a better behaved generalized option price function ${\cal
O}(x_j)$ and we introduce the weighted jump density
$\mu(x)=e^x\nu(x)$. Then the forward formula, expressing option
prices in terms of the model parameters, is given in the spectral
domain by ($\cal F$ denotes the Fourier transform)
\[{\cal FO}(v)=\frac{1-\exp(T(-\sigma^2(v-i)^2/2+i\gamma(v-i)+{\cal F}\mu(v)-\lambda))}{v(v-i)},\quad
v\in\R,
\]
cf. Equation (2.7) in \citeasnoun{Belomestny/Reiss:2006}. Our
observations are modeled as
\[ O_j={\cal O}(x_j)+\eps_j,\quad j=1,\ldots,n,\]
with i.i.d. and centered noise variables $(\eps_j)$. For the sake of
an easier presentation here, we assume that the real parameters
$(\sigma^2,\gamma,\lambda)$ are known such that the backward
formula, expressing the transformed jump density in terms of the
option prices, is given by
\[ {\cal F}\mu(v)=T^{-1}\log(1-v(v-i){\cal
FO}(v))+\sigma^2(v-i)^2/2-i\gamma(v-i)+\lambda,\quad v\in\R.
\]
We construct an empirical version $\tilde{\cal O}$ from the
observations $(O_j)$ (e.g. by linear interpolation) and obtain by
substitution in this formula an empirical version ${\cal
F}\tilde\mu(v)$ of ${\cal F}\mu(v)$, which satisfies for small noise
levels
\[\abs{{\cal F}(\tilde\mu-\mu)(v)}=T^{-1}\babs{\log\Big(\frac{1-v(v-i){\cal
F\tilde O}(v)}{1-v(v-i){\cal FO}(v)}\Big)}\approx
T^{-1}\babs{\frac{v(v-i){\cal F(\tilde O-O)}(v)}{1-v(v-i){\cal
FO}(v)}}.
\]
This first order analysis already reveals the ill-posedness of the
calibration problem: the higher the frequency $\abs{v}$, the more
the error ${\cal F(\tilde O-O)}(v)$ in the observation domain is
amplified by the factor $v(v-i)$ as well as by the denominator which
tends to zero for $\abs{v}\to\infty$. The nonlinearity is reflected
by a noise level which depends on the unknown true value ${\cal
FO}(v)$. A natural approach is to cut-off high frequencies and to
consider for $U>0$ the estimators
\[\widehat{\mu}_U(x)={\cal F}^{-1}({\cal F}\tilde\mu{\bf
1}_{[-U,U]})(x),\quad x\in\R.
\]

Let us write abstractly
\[ {\cal F}\tilde\mu(v)={\cal F}\mu(v)+\sigma(v)\xi(v),\quad v\in\R,\]
with $\sigma(v)$ denoting the noise level and $\xi(v)$ denoting the
normalized noise variable, that is the difference between empirical
and true value divided by its standard deviation, at frequency $v$.
Now we see the analogy with the sequence space model \eqref{EqModel}
analyzed before. In fact, the only difference is the continuity of
the spectral parameter $v$ instead of $k\in\N$ and the additional
difficulty is the much more complicated noise structure.
Nevertheless, the estimators $\widehat{\mu}_U$ are provably
rate-optimal for the right choice of the cut-off frequency $U$ (the
problem is severely ill-posed in general). In simulations the
standard data-driven choices of $U$ (e.g. Lepski's method, cross
validation) had a poor performance compared to the oracle choice,
mostly because of a very badly known noise structure. This is
exactly why the quasi-optimality criterion is of interest here.


\subsection{Experimental Setup}

In total we performed 1000 independent experiments. Each of them was
set up as follows (notation as in \citeasnoun{Belomestny/Reiss:2006}
where more details can be found):

\begin{figure}[t]
\begin{minipage}[b]{.48\linewidth}
\pgfdeclareimage[interpolate=false,width=\linewidth]{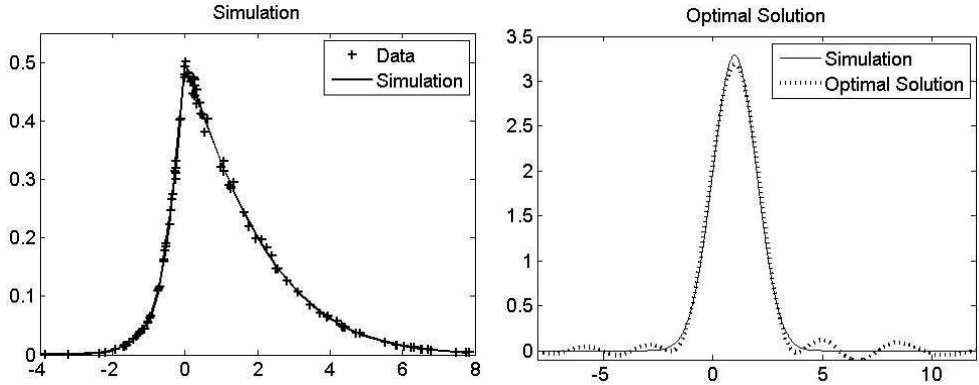}{data}
\centering\pgfuseimage{data}
\end{minipage}
\begin{minipage}[b]{.04\linewidth}
\end{minipage}
\begin{minipage}[b]{.48\linewidth}
\pgfdeclareimage[interpolate=false,width=\linewidth]{sol}{optimalsol}
\centering\pgfuseimage{sol}
\end{minipage}
\caption{left: Data, right: optimal solution}\label{Fig1}
\end{figure}

\begin{itemize}
\item $100$ design points $(x_i)$ were chosen at random, $50$ according to a standard normal,
$50$ according to a uniform distribution on $[-4,8]$.

\item Corresponding observations $(O_i)$ are generated by calculating the exact value and adding standard normal noise
of $3\%$ relative noise level, i.e. $O_i = {\cal O}(x_i)(1 + 0.03
e_i)$ with $e_i\stackrel{\text{iid}}{\sim} \mathcal{N}(0,1)$. The
other characteristics are chosen as follows (corresponding to the
so-called Merton model): volatility $\sigma=0.1$, jumps are standard
normal with intensity $\lambda=5$ (the jump density $\nu$ in Figure
\ref{Fig1} (right)), maturity $T=0.25$. The value of $\gamma$ is
obtained by imposing a martingale condition. In Figure \ref{Fig1}
(left) the true option price function ${\cal O}(\cdot)$ is shown
together with the observations as a function of $x$. Abstractly, one
can show that in this setting both, the bias and the variance term
have exponential decay respectively growth in $U$.

\item In total $60$ cut-off frequencies $(U_n)$ were chosen with step-width $0.8$, i.e. $U_n=0.8n$, $n=1,\ldots,60$.
\end{itemize}

\subsection{Parameter choice methods}

We will compare three different parameter choice methods. The
quasi-optimality and the Hardened Balancing Principle are treated in
this article and Lepski's method serves as a widely used benchmark.

\subsubsection*{Quasi-Optimality.}
We use exactly as in Definition \ref{DefParameter} with
$\chi(\cdot)=1$
\begin{equation*}
n_{qo} := \argmin_{n} \{ \| \hat x^{(n+1)} - \hat x^{(n)} \| \}.
\end{equation*}

\subsubsection*{Lepski-type method.}

Define
\begin{equation*}
f(n) := \max_{n<m\leq N} \left\{ 4^{-1} \|  \hat x^{(m)} - \hat
x^{(n)} \| / s(m) \right\},\quad n_{Lepski} := \argmin_n \{ \forall
m \ge n:\: f(m) \le \kappa \}.
\end{equation*}
For the choice of $\kappa$ one should theoretically use a quantity
depending on the noise level and larger than $1$. In practice,
however, it is observed that $\kappa \in [ 0.25, 0.75]$ gives
superior results and works in more or less any situation. Here we
choose $\kappa=0.75$.

\subsubsection*{Hardened Balancing Principle.}

We reuse the computed quantity $f(n)$ defined above and choose as regularization
parameter
\begin{equation*}
n_{HBP} := \argmin_{n} \{ f(n) \sqrt{s(n)} \}
\end{equation*}
Due to the definition of $f$ it can be seen as a stabilized version
of Definition \ref{DefParameter} with $\chi(n) = {s(n)}^{-1/2}$.
This stabilization mainly counters the effects of a subsampling with
inappropriately small spacing between the cut-off points.

\subsubsection*{Implementation details.}
There are certain fine points in the implementation, which are not
crucial, but yield slightly superior results. The important point
for the methods proposed here is the validity of Assumption
\ref{AssDecay}. Therefore we estimated $s(\cdot)$ out of $10$
independent data sets (each of them constructed as given above, i.e.
each of them with different design points and random error) and
removed the cut-off parameters $n\in \{ 17, 27, 36, 44, 46, 48 \}$
which seem to violate the assumption empirically.

The Lepski balancing principle and Hardened Balancing show undesirable properties when
$s(\max)/s(n) < 2$. Therefore we restricted the region of admissible regularization
parameters further to $n\in\{1, \ldots 43\}\setminus\{17, 27, 36\}$. The remaining
cut-off points were still used to calculate the values $f(n)$.

\subsection{Experiment}

\begin{figure}[t]
\pgfdeclareimage[interpolate=false,width=\linewidth]{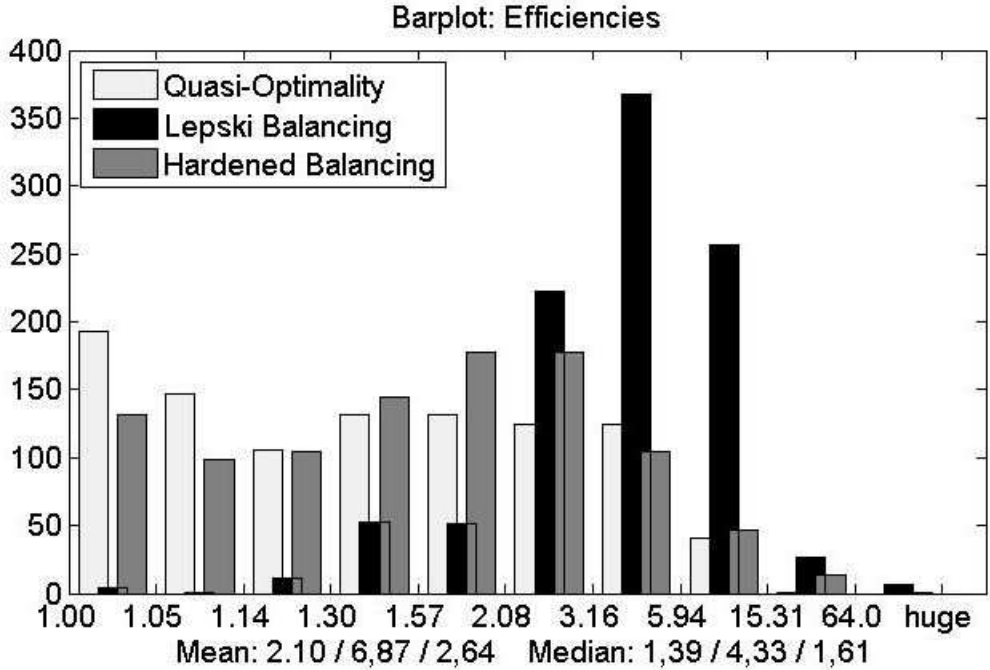}{efficiencies}
\pgfuseimage{efficiencies} \caption{Efficiencies}\label{fig2}
\end{figure}

In Figure \ref{fig2} the distribution of the efficiencies is
displayed, i.e. the ratio of errors between the unsupervised
solution selected by one of the methods and the optimal one, i.e.
\begin{equation*}
\operatorname{efficiency}(x,y,n_{chosen}) := \|x^{(n_{chosen})} -
x\| / \min_n \|x^{(n)} - x\|.
\end{equation*}
Note that the ratio cannot be smaller than $1$. Remark further that
it is much harder to obtain a low value for $\mathbb{E}
\operatorname{efficiency}(x,y,n_{chosen})$ compared to an oracle
criterion of the form $\mathbb{E}\|x^{(n_{chosen})} - x\| / \min_n
\mathbb{E} \|x^{(n)} - x\|$ because already one much better
estimator in the denominator can spoil the results in our case. The
bar plot of Figure \ref{fig2} bins the solutions according to their
efficiency. The bins are taken on a double exponential scale, ratios
larger than 64 were considered as huge.

As we can see, quasi-optimality and hardened balancing, as defined
here, are superior to the Lepski-balancing principle, which suffers
from a considerable number of bad results. Although in this
situation the hardened balancing principle performs slightly worse,
it has a particular advantage. Generally, it is relatively
insensitive towards a bad choice of the subsampling whereas already
one badly chosen distance in the subsampling could bring
quasi-optimality out of track. Yet, it should be recalled that it
relies significantly on an estimate of the stochastic noise level
$s(\cdot)$, obtained from additional data sets.

\section{Perspectives}\label{SectionPersp}

Numerical experiments indicate that the results can be carried over to other
regularization methods like Tikhonov regularization. Unfortunately, the correlation
structure induced by Tikhonov regularization makes an analysis in comparison to the
spectral cut-off regularization much more demanding and remains a topic of future
research.

Interestingly, a way to improve the numerical results for spectral
cut-off even further is to use a two-step procedure, where first a
regularization parameter is obtained from a very coarse subsampling
and then the region around the chosen regularization parameter is
reconsidered with a finer subsampling. These types of procedures
work provably well for change point detection problems (e.g.
\citeasnoun{Korostelev:1987}), a setting, which is to some extent
related with our search for the intersection point of $b(\cdot)$ and
$s(\cdot)$.

\section{Proofs}\label{SectionProofs}\label{SecProofProp}

In order to evaluate the probabilities of staying above or below the
threshold, we will make use of the following technical result.

\begin{lemma}\label{LemDev}
Let $Z=\sum_{k=1}^\infty\alpha_k^2\zeta_k^2$ with
$\sum_{k=1}^\infty\alpha_k^2=1$ and $\zeta_k\sim N(0,1)$ iid. Then
\[ \forall\,z\in (0,1):\;\PP(Z\le z)\le \exp\Big(\frac{1-z+\log(z)}{2\max_k\alpha_k^2} \Big),\qquad
\forall\,z>0:\;\PP(Z\ge z)\le \sqrt{2}e^{-z/4}.
\]
\end{lemma}

\begin{proof}
For all $\lambda>0$ and $z\in (0,1)$ we have
\begin{align*}
\PP(Z\le z) &=
\PP\Big(\exp\Big(-\lambda\Big(\sum_{k=1}^\infty\alpha_k^2\zeta_k^2\Big)\Big)\ge
e^{-\lambda
z}\Big)\\
&\le e^{\lambda
z }\prod_{k=1}^\infty\E[\exp(-\lambda\alpha_k^2\zeta_k^2)]\\
&=e^{\lambda z}\prod_{k=1}^\infty(1+2\lambda\alpha_k^2)^{-1/2}.
\end{align*}
Due to $(1+2x)^{-1/2}\le \exp(-\frac{\log(1+\eps)}{\eps}x)$ for
$x\in [1,1+\eps]$, $\eps>0$, we deduce for $\lambda=
\eps(2\max_k\alpha_k^2)^{-1}$ and $\eps=z^{-1}-1$:
\[ \PP(Z\le z)\le \exp\Big(\frac{\eps}{2\max_k\alpha_k^2}\Big(
z-\frac{\log(1+\eps)}{\eps}\Big)\Big)
=\exp\Big(\frac{1-z+\log(z)}{2\max_k\alpha_k^2} \Big).
\]
By Jensen's inequality we obtain for $z>0$
\begin{equation}
\PP(Z\ge z) \le
e^{-z/4}\E\Big[\exp\Big(\sum_{k=1}^\infty\alpha_k^2\zeta_k^2/4\Big)\Big]
\le e^{-z/4}\E[\exp(\zeta^2/4)]=\sqrt{2}e^{-z/4}.
\end{equation}
\end{proof}

\begin{proof}[of Lemma \ref{LemIneq}]
For each $n\ge 1$ and $\alpha>0$ we derive from the concavity of the
log-function and Jensen's inequality
\begin{align*}
&R_{\alpha}(\hat{x}^{(n)})^2=\tilde\E[\norm{\hat{x}^{(n)}-x}^{\alpha}]^{2/\alpha}\\
&\ge \exp\big(\tilde\E[\log(\norm{\hat{x}^{(n)}-x}^2)]\big)\\
&=\exp\Big(\tilde\E\Big[\log\Big(\sum_{k=1}^{\ell(n)}\tfrac{\sigma(k)^2\xi_k^2}{s(n)+b(n)}
+\sum_{k=\ell(n)+1}^\infty \tfrac{x_k^2}{s(n)+b(n)}\Big)\Big]+\log(s(n)+b(n))\Big)\\
&\ge
\exp\Big(\sum_{k=1}^{\ell(n)}\tfrac{\sigma(k)^2}{s(n)+b(n)}\tilde\E[\log(\xi_k^2)]
+\sum_{k=\ell(n)+1}^\infty \tfrac{\gamma(k)^2}{s(n)+b(n)}\tilde\E[\log(x_k^2/\gamma(k)^2)]\Big)(s(n)+b(n))\\
&=(s(n)+b(n))e^{\E[\log(\zeta^2)]}=e^{\E[\log(\zeta^2)]}R_2(\hat{x}^{(n)})^2,
\end{align*}
which gives the first inequality. For the second, we use
$w^{\alpha/2}/\Gamma(\frac{\alpha}{2}+1) \le e^w$ for $w\ge 0$,
Lemma \ref{LemDev} and $(1-w)^{-1/2}\le 2^w$ for $w\in (0,1/2]$ to
obtain for any $\alpha>0$ and $\lambda\in
(0,\min(\gamma(\ell(n)+1)^{-2},\sigma(\ell(n))^{-2})/4)$:
\begin{align*}
\tilde\E[\norm{\hat{x}^{(n)}-x}^{\alpha}] & =
\Gamma(\tfrac{\alpha}{2}+1) \lambda^{-\alpha/2}
\tilde\E[ \lambda^{\alpha/2}\norm{\hat{x}^{(n)}-x}^\alpha / \Gamma(\tfrac{\alpha}{2}+1)] \\
&\le \Gamma(\tfrac{\alpha}{2}+1)\lambda^{-\alpha/2}
\tilde\E[\exp(\lambda\norm{\hat{x}^{(n)}-x}^2)]\\
&=\Gamma(\tfrac{\alpha}{2}+1)\lambda^{-\alpha/2}\prod_{k=1}^{\ell(n)}
(1-2\lambda\sigma(k)^2)^{-1/2}\prod_{k=\ell(n)+1}^\infty
(1-2\lambda\gamma(k)^2)^{-1/2}\\
&\le
\Gamma(\tfrac{\alpha}{2}+1)\lambda^{-\alpha/2}\exp(2\log(2)\lambda(s(n)+b(n))).
\end{align*}
The choice $\lambda=(s(n)+b(n))^{-1}$ fulfills the above requirement
and thus yields
\begin{equation*}
R_{\alpha}(\hat{x}^{(n)})=\tilde\E[\norm{\hat{x}^{(n)}-x}^{\alpha}]^{1/\alpha}\le
(4 \Gamma(\tfrac{\alpha}{2}+1))^{1/\alpha}R_2(\hat{x}^{(n)}).
\end{equation*}
The last inequality follows from the two others together with the
fact that Assumption \ref{AssDecay} implies
\begin{equation*}
\min_nR_2(\hat{x}^{(n)})^2\ge \min(b(n^\#),s(n^\#+1))\ge
\tfrac{C_b^{-1}}{2} R_2(\hat{x}^{(n^\#)})^2.
\end{equation*}
\end{proof}

\begin{proof}[of Proposition \ref{PropLD}]
Let us first introduce the following quantities:
\begin{align*}
a(n)&:=\sqrt{\chi(n)}\sqrt{(s(n+1)-s(n))+(b(n)-b(n+1))},\\
D(n)&:=\chi(n)\norm{\hat{x}^{(n+1)}-\hat{x}^{(n)}}^2,\\
Z(n)&:=\sqrt{D(n)}/a(n).
\end{align*}
We shall only consider $n> n^\#$ because the case $n<n^\#$ is
symmetric (exchanging $s$ and $b$) and the case $n=n^\#$ is trivial.
For $m:=n-n^\#\ge 1$ we find by Assumptions \ref{AssDecay} and
\ref{AssDecay2}
\begin{align*}
\frac{a(n^\#)^2}{a(n^\#+m)^2}&\le \frac{\chi(n^\#)}{\chi(n^\#+m)}
\frac{(C_s-1)s(n^\#)+(C_b-1)b(n^\#+1)}{(c_s-1)s(n^\#+m)}\\
&\le {C_\chi}^m\; \frac{(C_s-1)s(n^\#)+(C_b-1)s(n^\#+1)}
{(c_s-1)c_s^{m}s(n^\#)}\\
&\le \frac{C_bC_s-1}{c_s-1} \left(\frac{c_s}{
{C_\chi}}\right)^{-m}=\rho(n).
\end{align*}
Using Lemma \ref{LemDev}, we infer for any $K>0$
\begin{align*}
\tilde\PP(n^\ast=n) &\le \tilde\PP(D(n)\le
D(n^\#))\\
&\le \tilde\PP(Z(n)^2 \le \rho(n) Z(n^\#)^2)\\
&\le \tilde\PP(Z(n)^2\le K\rho(n)\log(\rho(n)^{-1}))
+\tilde\PP(Z(n^\#)^2\ge
K\log(\rho(n)^{-1}))\\
& \le
\exp\left(\frac{1-K\rho(n)\log(\rho(n)^{-1})+\log(K\rho(n)\log(\rho(n)^{-1}))}
{2\max_k \frac{\sigma(k)^2 + \gamma(k)^2}{a(n)}} \right) \\&\qquad+ \sqrt{2} \exp(-K\log(\rho(n)^{-1})/4)\\
&\le
(Ke\rho(n)\log(\rho(n)^{-1}))^{a(n)^2/2\max_k(\sigma(k)^2+\gamma(k)^2)}
+
\sqrt{2}\rho(n)^{K/4}\\
&=\sqrt{2}\rho(n)^{K/4}+ (Ke\rho(n)\log(\rho(n)^{-1}))^{r(n)/2}.
\end{align*}
The choice $K=2r(n)$ yields the result.
\end{proof}

\begin{proof}[of Theorem \ref{MainTheorem}]

To prove that $n^\ast$ is well defined, we infer from the proof of
Proposition \ref{PropLD} that for all $n>n^\ast$
\[ \tilde\PP(D(n)\le
D(n^\#))\le \sqrt{2}\rho(n)^{r(n)/2}+
(2er(n)\rho(n)\log(\rho(n)^{-1}))^{r(n)/2}.
\]
This exponential decay implies
\[ \lim_{m\to\infty} \tilde\PP(\exists n\ge m:\: D(n)\le
D(n^\#))\le \lim_{m\to\infty}\sum_{n=m}^\infty \tilde\PP(D(n)\le
D(n^\#))=0,
\]
which means that the probability that the criterion $D(\cdot)$ is
larger for some $n\ge m$ than at $n^\#$ tends to zero as
$m\to\infty$, hence $n^\ast=\argmin_n D(n)$ is well defined with
probability one.

For the main assertion we use H\"older's inequality with
$p^{-1}+q^{-1}=1$, Lemma \ref{LemIneq} and Assumption \ref{AssDecay}
for any $\alpha>0$ to obtain:
\begin{align*}
&\tilde\E[\norm{\hat{x}^{(n^\ast)}-x}^{\alpha}]\\
&=\sum_{m=-n^\#+1}^\infty
\tilde\E[\norm{\hat{x}^{(n^\#+m)}-x}^{\alpha}{\bf
1}_{\{n^\ast=n^\#+m\}}]\\
&\le\sum_{m=-n^\#+1}^\infty
\tilde\E[\norm{\hat{x}^{(n^\#+m)}-x}^{\alpha p}]^{1/p}\tilde\E[{\bf
1}_{\{n^\ast=n^\#+m\}}^q]^{1/q}\\
&\le\sum_{m=-n^\#+1}^\infty
R_{\alpha p}(\hat{x}^{(n^\#+m)})^{\alpha}  \tilde\PP(n^\ast=n^\#+m)^{1/q}\\
&\le (4\Gamma(\tfrac{\alpha p}{2}+1))^{1/p} \sum_{m=-n^\#+1}^\infty
R_2(\hat{x}^{(n^\#+m)})^{\alpha} \tilde\PP(n^\ast=n^\#+m)^{1/q}\\
&\le (4\Gamma(\tfrac{\alpha
p}{2}+1))^{1/p}\big(2s(n^\#)\big)^{\alpha} \Big( \sum_{m=0}^\infty
C_s^{\alpha m/2} \tilde\PP(n^\ast=n^\#+m)^{1/q}+ \sum_{m=0}^{n^\#-1}
C_b^{\alpha m/2} \tilde\PP(n^\ast=n^\#-m)^{1/q}\Big).
\end{align*}
We choose $q:=(\frac{\alpha}{2}+\frac{r}{2}\min(\frac{\log(c_b
c_\chi)}{\log(C_b)}, \frac{\log(c_s /
C_\chi)}{\log(C_s)}))/\alpha>1$ and infer from Proposition
\ref{PropLD}
\begin{align*}
&\sum_{m=1}^\infty C_s^{\alpha m/2} \tilde\PP(n^\ast=n^\#+m)^{1/q}\\
&\le (\sqrt{2}+(2er)^{r/2})^{1/q}\left(\frac{C_bC_s-1}{c_s-1}
\log\left(\frac{c_s-1}{C_bC_s-1}
\frac{c_s}{C_\chi}\right)\right)^{r/2q}\sum_{m=1}^\infty C_s^{\alpha
m/2}\left(\frac{c_s}{C_\chi}\right)^{-rm/2q}m^{r/2q},
\end{align*}
which is finite by the choice of $q$ and the restriction on
$\alpha$. With a symmetric argument for $m\le 0$ we conclude that
$\tilde\E[\norm{\hat{x}^{(n^\ast)}-x}^{\alpha}]$ is bounded by a
multiple of $s(n^\#)^{\alpha/2}$, which has the order of
$R_2(\hat{x}^{(n^\#)})^{\alpha}$. We eventually obtain with
constants $K=K(c_s,c_b,C_s,C_b,r,\alpha)$,
$K'=K'(c_s,c_b,C_s,C_b,r,\alpha)$ (note that $p$ and $q$ depend on
$\alpha$ and the remaining constants):
\begin{align*}
R_{\alpha}(\hat{x}^{(n^\ast)}) &\le K' R_2(\hat{x}^{(n^\#)})\le
K\min_{n\ge 1}R_{\alpha}(\hat{x}^{(n)}).
\end{align*}
\end{proof}

\section*{Acknowledgements}

The first author gratefully acknowledges the financial support by
the Upper Austrian Technology and Research Promotion. We are
grateful for the constructive criticism by two anonymous referees.


\bibliographystyle{economet}
\bibliography{bibliography}

\end{document}